\documentclass[11pt]{article}
\usepackage{amsmath}
\usepackage{amssymb}
\usepackage{amsfonts}
\usepackage{eucal}
\usepackage[usenames]{color}
\usepackage{graphicx}
\numberwithin{equation}{section}
\newtheorem{Theorem}{Theorem}[section]
\newtheorem{Proposition}[Theorem]{Proposition}
\newtheorem{Definition}[Theorem]{Definition}

\newtheorem{Remark}[Theorem]{Remark}

\begin{document}
\title{One-dimensional reflected diffusions
with two boundaries and an inverse first-hitting problem}
\author{Mario Abundo\thanks{Dipartimento di Matematica, Universit\`a  ``Tor Vergata'', via della Ricerca Scientifica, I-00133 Rome, Italy.
E-mail: \tt{abundo@mat.uniroma2.it}}
}
\date{}
\maketitle

\begin{abstract}
\noindent We study an inverse first-hitting problem for a one-dimensional, time-homogeneous diffusion $X(t)$ reflected between two boundaries $a$ and $b,$
which starts
from a random position $\eta.$ Let  $a \le S \le b$ be  a given threshold, such that $P( \eta \in [a,S])=1,$ and
$F$ an assigned distribution function.
The  problem consists of finding the distribution of $\eta$ such that the first-hitting time of $X$  to $S$ has distribution $F.$
This is a generalization of the analogous problem for ordinary diffusions, i.e. without reflecting, previously considered by the author.

\end{abstract}


\section{Introduction}
Reflected diffusion processes with one or two  boundaries play an important role in a
variety of applications ranging from Economics, Finance,
Queueing, and  Mathematical Biology.
As far as Economics and Finance are concerned, we mention e.g. certain models for the currency exchange rate
(\cite{ball:amf98}, \cite{bert:amer92}, \cite{dejong:appe94}, \cite{krugman:quart91}, \cite{sven:jmon91});
for other applications in Economics, and Insurance, see e.g. \cite{qin:2012}, \cite{veer:comec04}. As for Queueing theory,
diffusions with reflecting boundaries arise
as heavy-traffic approximations of queueing systems (see e.g. \cite{abate:aap87a}, \cite{abate:aap87b}, \cite{harr:85}
for reflected Brownian motion, and \cite{sri:96}, \cite{ward:que03a}, \cite{ward:que03b} for reflected
Ornstein-Uhlenbeck (OU) process).
Reflected OU process appears also  in certain models from Mathematical Biology (see e.g. \cite{ricsac:87}).
For further applications, see e.g.
\cite{linet:2005} and references therein.
In all these, the knowledge of the distribution of the
first-passage-time (FPT) of the reflected diffusion through an
assigned  barrier is very important, in order to obtain a more
precise insight of the modeled phenomenon. Although FPT problems
have been studied mostly for ordinary diffusions, i.e. without
reflecting (see e.g. \cite{abundo:pms00}, \cite{abundo:pms97},
\cite{abundo:mcap08}, \cite{abundo:saa03}, \cite{abundo:stapro02},
\cite{abundo:sta02}, \cite{daniel:jap69}, \cite{darling:ams53},
\cite{durbin:jap71}, \cite{martin:jap98}, \cite{ric:lam90}, \cite{salm:aap88}), more
recently (see e.g.  \cite{chuan:2012}, \cite{lijun:2006},
\cite{qin:2012} ) some results appeared about the FPT of a
one-dimensional reflected  diffusion, through a threshold $S.$
\par In this paper, we  focus on FPT problems for a
one-dimensional, temporally homogeneous reflected diffusion
process $X(t)$ with boundaries $a$ and $b,$  which is the solution
of the stochastic differential equation with reflecting boundaries
(SDER):
\begin{equation} \label{eqdiffur}
\begin{cases}
dX(t)= \mu (X(t)) dt + \sigma (X(t)) d B_t  + dL_t - dU_t \\  X(0) = \eta \in [a,b]
\end{cases} \ ,
\end{equation}
where
$B_t$ is standard Brownian motion,
the initial position $\eta$ is a random variable, independent of  $B_t,$
 $L= \{ L_t \}$ and $U= \{ U_t \},  \ t \ge 0,$
are the {\it regulators} of points $a$ and $b,$ respectively, namely the local times of $X$ at $a$ and $b.$
The processes $L$ and $U$ are uniquely determined by the following properties (see e.g. \cite{harr:85}): \par
(i) both  $L_t$ and $U_t$ are continuous nondecreasing processes with $L_0=U_0=0;$ \par
(ii) $X(t) \in [a,b]$ for every time $t \ge 0;$ \par
(iii) $L$ and $U$ increase only when $X=a$ and $X=b,$ respectively, that is,
for $t\ge 0,$ \par
$ \ \ \ \  \int _0 ^t {\bf 1 } _ { \{ X(s)=a \} } dL_s =L_t$ and $\int _0 ^t {\bf 1 } _ { \{ X(s)=b \} } dU_s =U_t.$ \bigskip

Under certain mild regularity conditions on the coefficients $\mu( \cdot)$ and $\sigma (\cdot)$ (see e.g. \cite{lions:1984}),
for fixed initial value, the SDER \eqref{eqdiffur}
has a unique strong solution $X(t)$ which remains in the interval $[a,b]$ for every time $t \ge 0.$ For this reason,
$X(t)$ is also called a {\it regulated diffusion } between $a$ and $b$.
\par
When $X$ is Brownian motion ($\mu \equiv 0, \ \sigma \equiv 1 )$ with only reflecting boundary $a=0,$
the SDER \eqref{eqdiffur} in its integral form becomes the Skorohod's
equation for  reflected
BM, that is,  $X(t) = X(0) + B_t + l_t \ ,  X(0) \ge 0,$
where $X(t)$ and $l_t$ have to be found as continuous functions under the conditions that $X(t) \ge 0, \ l_0=0$ and $l_t$ is non-negative, nondecreasing,
flat outside
the set $\{ t: X(t) =0 \}.$
The Skorohod's equation has a unique solution
for all but a negligeable class of Brownian paths, and it is given by (see \cite{mckean:63}, \cite{mckean:69},
\cite{sko:62}):
\begin{equation}
X(t)=
\begin{cases}
X(0)+ B_t \ , &  0 \le t \le T   \\ X(0) + B_t + l_t  \ , & t > T
\end{cases} \ ,
\end{equation}
where $l_t = \sup _ { s \in [0,t]} B_s ^- $ and $T = \inf \{ t >0: B_t \le 0 \}$ (here $B_t^- = - \min \{ B_t, 0 \} ).$
The process $X$ is identical in law to the diffusion with generator $ \widetilde {\cal L} \phi  = \frac 1 2 \phi '' $
subject to the reflecting boundary condition
$ \phi '_+ (0)= \lim _{ t \rightarrow 0 ^+} \frac 1 t \ (\phi (t)- \phi (0) )  =0 .$
Moreover (see e.g. \cite{ikwa:sde81}, \cite{mckean:63}),  $l_t$ is equivalent in law to
$ L_t= \lim _{ \epsilon \rightarrow 0 ^+} \frac 1 { 2 \epsilon} \int _ 0 ^t {\bf 1 } _ { [0, \epsilon)} ( X (s)) ds , $ i.e.
the local time of $X$ at zero.
\bigskip

Let $X$ be the solution of the SDER \eqref{eqdiffur}; if $S \in [a,b]$ is a threshold such that  $P( a \le \eta
\le S ) =1,$ we consider the FPT of $X$ through $S:$
\begin{equation} \label{FPT}
\tau _S  = \inf \{t >0 : X(t) = S \}
\end{equation}
and we denote by
$ \tau _S (x) = \inf \{t >0 : X(t) = S  | \eta =x\} $
the FPT of $X$ through $S$ with the condition that $\eta =x.$
We assume that  $\forall x \in [a,S], \ \tau_S (x)$ is finite with probability one and
that it possesses a density $f(t|x).$ \par
Usually, in the FPT problem the initial position $X(0)$ is assumed to be deterministic and fixed to a value $x,$
then the  {\it direct}  problem
consists of finding the FPT density $f(t|x)= f_S (t|x)$ or
the moments of   $\tau _S (x).$ Notice that, while for ordinary diffusions a certain number of results about the direct FPT problem is available
in the literature (see the papers cited above), for reflected diffusion processes only  a few papers appeared on this topics
(see e.g. \cite{chuan:2012}, \cite{lijun:2006}, \cite{qin:2012}). \par\noindent
The {\it inverse} FPT problem for diffusions generally focuses on determining
the barrier  $S,$ when  $f(t|x)$ is given (see e.g.
\cite{abundo:saa06}, \cite{zuc:aap09}).
Since  we assume that the initial position $\eta $ is random, we consider a  slight modification of the problem,
that is the following inverse first-passage-time  (IFPT) problem. \par\noindent
Let be $S \in [a,b]$ an assigned barrier, and suppose that the initial position  $X(0) = \eta $ is independent of $B_t$
and  $P( a \le \eta \le S  )=1 ;$
for a given distribution $F,$ our aim is: \bigskip

 {\it to find the density g of $ \eta $ (if it  exists) for  which
it  results  $P(\tau_S  \le t) = F(t)$}. \bigskip

\noindent This IFPT problem has interesting applications in
Mathematical Finance, in particular in credit risk modeling, where
the FPT represents a default event of an obligor (see e.g.
\cite{jackson:stapro09}),  in Biology, specially in the
framework of diffusion models for neural activity (see e.g.
\cite{lansky:mbs89}), and in Queueing theory (see e.g. \cite{abate:aap87a}, \cite{abate:aap87b}, \cite{harr:85}).
For ordinary diffusions, it was studied in
\cite{jackson:stapro09} in the  case of Brownian motion, while
some extensions to more general processes were obtained in
\cite{abundo:stapro13}, \cite{abundo:stapro12}.
\par The paper is organized as  follows:
Section 2 contains the formulation of the problem and the main
results, in Section 3 some explicit examples are reported.

\section{Notations, formulation of the problem and main results}
For $a < b ,$ let  $X(t)$ be  a one-dimensional, time-homogeneous diffusion with reflecting boundaries $a$ and $b,$
which is solution of the SDER (1.1), whose drift $\mu(\cdot)$ and infinitesimal variance $\sigma ^2(\cdot)$
are supposed to be sufficiently regular functions (see \cite{lions:1984}), in order to guarantee the existence and uniqueness of the strong
solution, for fixed initial condition. Then, $X(t)$ turns out to be a time-homogeneous strong Markov process with infinitesimal generator:
\begin{equation}
{\cal L} \phi (x) = \mu(x) \phi '(x) + \frac 1 2 \sigma^2 (x) \phi '' (x), \ x \in (a,b)
\end{equation}
acting on $C^2-$functions $\phi$ on $(a,b).$
\par
Let $S \in [a, b]$  a given barrier;  if $x \in [a,S],  $ we denote by
$\tau _ { x \uparrow S } (x) = \inf \{ t >0: X(t) = S | X(0)=x \}$
the first-hitting time of $X$ to $S,$ with the condition $X(0)=x,$
namely the FPT of $X$ through $S,$ ``from below''.
In analogous manner, if $x \in [S,b], $ we denote by
${\tau}  _{ x \downarrow S } (x) = \inf \{ t >0: X(t) = S | X(0)=x \}$
the FPT of $X$ through $S,$ ``from above''. To simplify notations, from now on, we will denote $\tau _ { x \uparrow S } (x) $ by $\tau_S (x)$ and
$\tau _ { x \downarrow S } (x) $ by $\widetilde {\tau}  _S (x).$
We will suppose that $\tau _S(x)$ and $\widetilde {\tau}  _S (x)$
are a.s. finite, for every fixed $x .$ \par\noindent
Preliminarily, we recall some facts concerning the direct first-hitting problem for $X.$

\subsection{The Laplace transform and the moments of the first-hitting time of $X$ to $S$}
The following result holds:
\begin{Theorem} {\rm (\cite{chuan:2012})} \
Let $X$ be the solution of the SDER \eqref{eqdiffur} with deterministic and fixed initial condition $X(0)=x,$ and
let $S \in [a,b].$ For  $x \in [a,S]$ and
$ \theta \ge 0,$ suppose that $u(x)= u_ \theta (x)$ satisfies the following equation:
\begin{equation} \label{laplacetaubeloweq}
\begin{cases}
{\cal L} u (x) = \theta u(x), \ x \in (a,S) \\
u'(a) =0
\end{cases} \ .
\end{equation}

Then, if $u(S) \neq 0$ for $S \in [x,b],$
the Laplace transform of $\tau_S(x)$ is explicitly given by:
\begin{equation} \label{laplacetaubelow}
E \left (e ^ {- \theta \tau_S(x) } \right ) = \frac {u(x) } {u(S)} \ .
\end{equation}
\par
In analogous manner, for $x \in [S,b]$ and $\theta \ge 0,$  let $v(x)= v_ \theta (x)$ be the solution of the problem:
\begin{equation} \label{laplacetauaboveeq}
\begin{cases}
{\cal L } v (x) = \theta v(x), \ x \in (S,b) \\
v'(b) =0
\end{cases} \ .
\end{equation}

Then, if $v(S) \neq 0 $ for $S \in [a,x],$
the Laplace transform of $\widetilde{ \tau} _S(x)$ is explicitly given by:
\begin{equation} \label{laplacetauabove}
E \left (e ^ {- \theta \widetilde {\tau }_S(x) } \right ) = \frac {v(x) } {v(S)} \ .
\end{equation}

\end{Theorem}
\par
\hfill $\Box$

\begin{Remark} {\rm As already noted in \cite{chuan:2012}, though neither $u$ and $v$  is unique, each of their
ratios is unique.
Thus, the Laplace transform of the first-hitting time of $X$ to $S$ can be obtained by solving a differential problem of the
second order. }
\end{Remark}

Let us consider e.g. $\tau_S(x),$ namely the FPT from below with the condition  $X(0)=x \in [a,S] , $ and denote by $M(\theta, x ) = u(x)/u(S)$ the Laplace transform
of $\tau_S(x)$
and by $T_n(x) =E[(\tau_S(x))^n]$ its moment of order $n, \ (n =1, 2, \dots) ;$ as well-known,
$T_n (x)= (-1)^n \frac { \partial ^n  } { \partial \theta ^n} M(\theta, x) \big | _ { \theta =0 } \ ,$ if it exists finite.
For fixed $\theta \ge 0,$ one has $u(x)= u(S) M(\theta, x);$ so the problem \eqref{laplacetaubeloweq} can be written as:
\begin{equation} \label{equationforM}
\begin{cases}
{\cal L} M(\theta , x) = \theta M( \theta , x), \ x \in (a,S)  \\ \frac {\partial } {\partial x } M ( \theta ,x ) \big |_ {x =a} =0
\end{cases} \ ,
\end{equation}
where the operator ${\cal L }$ acts on $M$ only as a function of $x.$

Therefore, by taking the $n-th$ derivative of $M ( \theta, x)$ with respect to $\theta$
in both members of equation \eqref{equationforM}, and calculating it at
$\theta =0,$ we obtain:

\begin{Proposition}
For $n=1, 2, \dots ,$ the $n-th$ order moments of $\tau_S(x),$
if they exist finite, are the solutions to the problems:
\begin{equation} \label{momentstaubelow}
\begin{cases}
{\cal L} T_n(x) = -n T_{n-1} (x), \ x \in (a,S) \\ T_ n (S) =0, \ T '_n (a) =0
\end{cases} \ ,
\end{equation}
where $T_0(x) \equiv 1.$
\end{Proposition}
\par
\hfill $\Box$
\bigskip

\noindent Analogous equations hold for the $n-th$ order moments of $\widetilde {\tau} _S(x),$
say $\widetilde T _n(x),$
 but the equations hold for $x \in (S,b)$ and the second boundary conditions have to be
replaced with $\widetilde T '_n(b)=0.$
\bigskip
\newline

Now, we derive the explicit solutions of problems  \eqref{laplacetaubeloweq} and  \eqref{momentstaubelow}
for the Laplace transform and the moments of  $\tau _S(x),$
in the case of reflected
 Brownian motion with drift $\mu,$ that is
the diffusion $X^ {(\mu )} $ with reflecting boundaries $a$ and $b,$ having infinitesimal generator
${\cal L} ^ {( \mu )}  \phi (x)= \mu \phi '(x) + \frac 1 2 \phi '' (x) .$
For $x \in [a,S],$ denote by
$\tau_S ^ {( \mu ) } (x)$ the FPT of $X^ {(\mu )} $ through $S$ from below, and by $T_n ^ { (\mu ) } (x)$ its moment of order $n;$
solving \eqref{laplacetaubeloweq} by quadratures and using  \eqref{laplacetaubelow}, we get that the Laplace transform of  $\tau_S ^ {( \mu )} (x)$ is:
$$
M^ {( \mu)} ( \theta , x) := E \left ( e ^ {- \theta \tau_S ^ {( \mu) } (x) } \right )
$$
\begin{equation} \label{laplacetaubelowexplicit}
= e^{ - (S-x)( \sqrt { \mu ^2 + 2 \theta } - \mu )  } \cdot
 \frac {\theta e^{-2 (x-a)\sqrt {\mu ^2 +2 \theta } } + \mu ^2 + \theta + \mu \sqrt { \mu ^2 + 2 \theta } }
{ \theta e^{-2 (S-a) \sqrt {\mu ^2 +2 \theta } } + \mu ^2 + \theta + \mu \sqrt { \mu ^2 + 2 \theta } } \ .
\end{equation}
For $a \rightarrow - \infty $ the right-hand member of \eqref{laplacetaubelowexplicit} tends to
$e^{-(S-x) ( \sqrt {\mu ^ 2 + 2 \theta} - \mu )},$ which is the well-known expression of the Laplace transform of the first-hitting time of
ordinary Brownian motion with drift $\mu $ to $S,$ when starting from $x < S.$ \par
Taking the limit as $\mu$ goes to zero in \eqref{laplacetaubelowexplicit}, we obtain:
\begin{equation} \label{laplacetaubelowmu0}
M^ {( 0)} (\theta , x) = E \left ( e ^ {- \theta \tau_S ^ {( 0)} (x) } \right )=
\frac {e^{ -x \sqrt {2 \theta} } + e ^{- (2a -x) \sqrt {2 \theta } } }
{e^{ -S \sqrt {2 \theta} } + e ^{- (2a -S) \sqrt {2 \theta } } } \ .
\end{equation}
In the special case $a=0,$ the expression above writes:
\begin{equation}
 \frac {e^ { -x \sqrt {2 \theta} } + e^ { x \sqrt {2 \theta} }  } {e^ { -S \sqrt {2 \theta} } + e^ { S \sqrt {2 \theta} }  } =
\frac {cosh (x \sqrt { 2 \theta} ) } {cosh (S \sqrt { 2 \theta} ) } , \ x \in [0,S].
\end{equation}
Then, Laplace transform inversion yields that, for $a=0$ and $x \in [0,S],$ the density of $\tau_S ^ {( 0)} (x) $  is
(cf. e.g. \cite{darling:ams53}, \cite{qin:2012}):
\begin{equation}
f^ {(0)} (t| x  )= \frac \pi {S ^2 } \sum _ {
k=0} ^ \infty (-1)^k \left (k + \frac 1 2 \right ) \cos \left [
\left (k + \frac 1 2 \right ) \frac { \pi x } S \right ] \exp
\left [ -  \left (k + \frac 1 2 \right ) ^2 \frac {\pi ^2 t } {2
S ^2 } \right ], \ t\ge 0 .
\end{equation}
\bigskip

\noindent Solving \eqref{momentstaubelow} with $n=1 $ by quadratures, we obtain the mean of $\tau_ S ^ {( \mu) } (x):$
\begin{equation} \label{T1}
T_1 ^ {( \mu) } (x)= \frac 1 { 2 \mu ^2 } \left [ e^ { 2 \mu (a-S)} - e^ { 2 \mu (a-x) } \right ] + \frac {S-x } { \mu}, \ x \in [a,S].
\end{equation}
Letting  $\mu$ go to zero, we obtain:
\begin{equation} \label{T10}
T_1 ^ {( 0 )} (x) = - x^2 +2ax +S (S -2a), \ x \in [a,S].
\end{equation}
As for the second order moment, by solving \eqref{momentstaubelow} with $n=2 $ by quadratures, we get:
\begin{equation}
T_2 ^ {( \mu ) } (x)=  \frac {x^2 } { \mu ^2 } - \frac x { \mu ^3 } \left ( e^ { 2 \mu (a-S) } +1 + 2S \mu + e ^ { 2 \mu (a-x) } \right ) +
c_1 + c_2 e^ { -2 \mu x },
\end{equation}
where
$$ c_2 = \frac {e^{2 \mu a } } {2 \mu ^4 } \ \left [ 4 a \mu - e^{2 \mu (a-S)} -2 -2 S \mu \right ],$$
and
$$ c_1 = - c_2 e^{- 2 \mu S} + \frac {S } {\mu ^3 } \ \left [2 e^{2 \mu (a-S)} +1 + S \mu \right ] .   $$
For $\mu =0 ,$  we obtain:
$$ T_2 ^ {( 0 ) } (x)= \frac {x^4} 3 - \frac 4 3 a x^3 -2S(S-2a)x^2 + Ax + B , \ x \in [a,S],$$
where
$$A= \frac 8 3 a^3 + 4aS(S-2a), \ B= \frac 5 3 S^4 - \frac {20} 3 a S^3 + 8 S^2 a^2 - \frac 8 3 a^3 S  .$$
In particular, for $a=0,$  we get
$ T_1 ^ {( 0 )} (x)= -x^2 +S^2 , \ T_2 ^ {( 0 )} (x)= \frac {x^4} 3 -2S^2x^2 + \frac 5 3 S^4 ,$
and so the variance of $\tau ^{(0)} _S (x)$ is $Var \left ( \tau ^{(0)} _S (x) \right )= \frac 2 3 \left (S^4-x^4 \right ), \ x \in [0,S].$ \bigskip

Explicit formulae for the Laplace transform  of the first-hitting time to a barrier $S$ are  known also
for reflected OU process,  reflected Bessel process and some other processes (see \cite{chuan:2012}, \cite{lijun:2006}),
but they involve special functions. An explicit spectral representation of the hitting time density was found in
\cite{qin:2012} for reflected BM, and in \cite{linet:2004a}, \cite{linet:2004b} for Cox-Ingersoll-Ross (CIR) and OU processes.

\subsection{The IFPT problem for reflected Brownian motion with drift}
For a given barrier $S \in [a,b],$ and $X(0)= \eta \in [a,S],$
let $\tau_S$ be the
FPT of $X(t)$  through $S$ from below (see \eqref{FPT}), and  suppose that $\tau _S(x)$ (i.e. the
FPT  conditional to $\eta =  x )$ is a.s.
finite for every $x \in [a,S] ,$ and it possesses a density $f (t |x).$
Moreover, we suppose that the initial position $\eta$ has a
density $g(x)$ with support  $(a,S);$  for $\theta \ge 0$ we denote by
$\widehat f(\theta | x)= \int _ 0 ^{ + \infty } e ^{ - \theta x} f(t|x) dt$ the Laplace transform of $f( t |x)$
and by
$\widehat g ( \theta )= \int _a^S e^ { - \theta x } g(x)
dx $ the (possibly bilateral) Laplace transform of $g.$
Then,  the density
of $\tau_S$ is obtained as $ f(t) = \int _ a ^ S f(t |x ) g(x) dx $
and taking the Laplace transform on both sides we get:
\begin{equation} \label{laplacegingeneral}
\widehat f(\theta)= \int _ a ^ S  \widehat f (\theta |x ) g(x) dx \ .
\end{equation}

Now, we go to solve the IFPT problem, in the case when $X= X ^ {(\mu)}$ is reflected BM with drift $\mu ,$ between
the boundaries $a$ and $b.$
\par\noindent
For a given FPT distribution function $F$ (or equivalently for a given FPT density $f=F')$
our aim is to find the density $g$ of the random initial position $\eta,$ if it exists, such that
$P( \tau_S  \le t) = F(t) .$ We are able to obtain the following:

\begin{Theorem} \label{primoteorema}
For $S \in [a,b],$  let $X^ {(\mu)}$ be BM with drift $\mu ,$ reflected between the boundaries $a$ and $b$ and
starting from the random position $\eta \in [a, S];$ suppose that the FPT of $X^ {(\mu)}$ through $S$ from
below has an assigned probability density  $f$ and denote by
$ \widehat f (\theta ) = \int _0 ^ \infty e ^ {-\theta t} f(t) dt , \ \theta \ge 0 ,$
the Laplace transform of $f.$ Then, if
there exists a solution $g$ to the IFPT problem for $X^ {(\mu)}$, its Laplace transform
$ \widehat g $ must satisfy the equation:
$$
 \widehat f (\theta) = \left [ \theta e^{ 2a \sqrt {\mu ^2 + 2 \theta } } \widehat g  \left (\sqrt {\mu ^2 + 2 \theta } + \mu \ \right ) +
 \left ( \mu ^2 + \theta + \mu \sqrt { \mu ^2 + 2 \theta } \ \right ) \widehat g   \left (\mu - \sqrt {\mu ^2 + 2 \theta } \ \right ) \right ] \times
$$
\begin{equation} \label{laplaceg}
 \left [ \theta e^{ -S ( \sqrt { \mu ^2 + 2 \theta } - \mu ) } e^{ 2a \sqrt { \mu ^2 + 2 \theta } } +
\left ( \mu ^2 + \theta + \mu \sqrt { \mu ^2 + 2 \theta } \ \right ) e^{ S( \sqrt { \mu ^2 + 2 \theta } + \mu ) } \right ] ^{-1}, \ \theta \ge 0 .
\end{equation}
In particular, if $\mu =0,$ the above formula takes the form:
\begin{equation} \label{laplacegmu0}
\widehat f (\theta) = \frac {\widehat g  ( \sqrt {2 \theta } ) + \widehat g  ( - \sqrt {2 \theta } ) e^{ -2a \sqrt {2 \theta } } }
{e^{ -S \sqrt {2 \theta} } + e ^{ (S-2a) \sqrt {2 \theta } }  } , \ \theta \ge 0 .
\end{equation}
Furthermore, if $\mu =0$ and if the density $g $ is required to be symmetric with respect to $(a+S)/2,$
then:
\begin{equation} \label{laplacegmu0symm}
\widehat g  ( \theta ) = \frac {e^{-S \theta } + e ^{ -(2a - S) \theta } } {1 + e^{ (S-a) \theta } } \
\widehat f \left ( \frac {\theta ^2 } 2 \right ) , \ \theta \ge 0 .
\end{equation}
\end{Theorem}
{\it Proof.} By using \eqref{laplacetaubelowexplicit} with $M^ {(\mu)} (\theta, x) = \widehat f (\theta |x),$ we have:
$$ \widehat f ( \theta)= \int _ a ^S \widehat f (\theta |x ) g  (x) dx $$
$$= e^{ -S ( \sqrt { \mu ^2 + 2 \theta } - \mu ) } \left [ \theta e ^{ -2(S-a) \sqrt { \mu ^2 + 2 \theta }} + \mu ^2 + \theta +
\mu \sqrt { \mu ^2 + 2 \theta } \ \right ] ^{-1} \times $$
$$ \int _a ^S e^{x ( \sqrt { \mu ^2 + 2 \theta } - \mu )} \left [\theta e ^{ -2(x-a) \sqrt { \mu ^2 + 2 \theta }} + \mu ^2 + \theta +
\mu \sqrt { \mu ^2 + 2 \theta } \ \right ] g(x)  \ dx  \ .$$
The integral can be written as
$$ \int _ a ^S \left [ \theta e^{ -x( \sqrt {\mu ^2 + 2 \theta } + \mu ) } e^{ 2a \sqrt {\mu ^2 + 2 \theta }} +
e^{x (\sqrt {\mu ^2 + 2 \theta } - \mu ) }  \left (\mu ^2 + \theta + \mu \sqrt {\mu ^2 + 2 \theta } \ \right )  \right ] g  (x) dx $$
$$= \theta e^{ 2a \sqrt {\mu ^2 + 2 \theta }} \widehat g ( \sqrt {\mu ^2 + 2 \theta } + \mu ) +
\left (\mu ^2 + \theta + \mu \sqrt {\mu ^2 + 2 \theta } \ \right ) \widehat g  ( \sqrt {\mu ^2 + 2 \theta } - \mu ).$$
Thus, by inserting this in the formula above,  \eqref{laplaceg} follows, after some manipulation.
Formula \eqref{laplacegmu0} is soon obtained, by taking $\mu =0;$ moreover, if one seeks
that the density $g $ is symmetric with respect to $(a+S)/2,$ namely
$\widehat g ( - \theta) = e^{(a+S) \theta } \widehat g  ( \theta ),$
one can explicit formula \eqref{laplacegmu0}
with respect to the Laplace transform of $g,$ and \eqref{laplacegmu0symm} follows.
\par
\hfill $\Box$

\begin{Remark}
{ \rm
Notice that, as in the case of the IFPT problem for ordinary diffusions, i.e. without reflecting (see e.g.
\cite{abundo:stapro13}), \cite{abundo:stapro12}, \cite{jackson:stapro09} )
the function $\widehat g$ may not be the Laplace transform of some
probability density function; in that case the IFPT problem has no solution.
This is the reason why Theorem \ref{primoteorema}
is formulated in a conditional form.
This kind of difficulty in showing the existence of a solution to an inverse first-passage time  problem is common to
another type of inverse
problem (see e.g. \cite{abundo:saa06}), in which one has to find the shape of the moving barrier $S(t)$ in such a way that the
FPT  of $X(t)$ over $S(t)$ for deterministic fixed initial condition, has a predetermined distribution $F.$
A numerical solution to that inverse problem was found in \cite{abundo:saa06}; the existence of the
solution is, at the moment, a still open problem  (see e.g. \cite{zuc:aap09}).
\par\noindent
If we replace reflected drifted Brownian motion with a more general reflected diffusion,
we still obtain \eqref{laplacegingeneral}, of course; however, even if
the explicit form of the Laplace transform $\widehat f (\theta |x)$ of the conditional FPT density $f(t |x )$ is available, in general
the expression \eqref{laplacegingeneral} cannot be put in
relation with the Laplace transform of $g$ calculated in some point (this happens for  reflected drifted Brownian motion,
thanks to the particular form of the Laplace transform of the FPT,  which depends on $x$ only by means of
exponentials of linear functions of $x$ (see \eqref{laplacetaubelowexplicit}).
}
\end{Remark}

If $\widehat f (\theta)$ is analytic in a neighbor of $\theta= 0,$ then the $k-$th order moments of $\tau_S$ exist
finite and they are obtained in terms
of $\widehat f (\theta)$
by $E(\tau_S ^k) = (-1)^k \frac {\partial ^k } {\partial \theta ^k } \widehat f(\theta)|_{\theta =0}.$
The same thing holds for the moments of $\eta,$
if $\widehat g (\theta)$ is analytic.
In the case of reflected drifted Brownian motion, by inserting the expression of $E( \tau ^{(\mu )} _S (x))$ given by \eqref{T1} into the equation
$$ E( \tau ^{ (\mu )}_S ) = \int _a ^S E( \tau ^{(\mu)}_S (x)) g(x) dx ,$$
after some manipulation, we obtain:
\begin{equation} \label{meantaumu}
E ( \tau ^{(\mu )} _S) =  \frac 1 \mu E( S - \eta ) - \frac 1 { 2 \mu ^2 } \ e^{ 2 \mu a } E \left (e^{-2 \mu  \eta } -e^{ -2 \mu S }  \right ) .
\end{equation}
Since it must be $E ( \tau ^{(\mu )} _S) \ge 0,$
we obtain the compatibility condition:
\begin{equation} \label{eqcompatibility1}
\frac 1 \mu E( S - \eta ) \ge \frac 1 { 2 \mu ^2 } e^{ 2 \mu a } E \left ( e^{ - 2 \mu \eta } - e ^{ -2 \mu S } \right ) .
\end{equation}
In particular, if $\mu > 0$ we have
$-\frac 1 { 2 \mu ^2 } e^{ 2 \mu a} E \left ( e^{ - 2 \mu \eta } - e ^{ -2 \mu S } \right ) \le 0,$ because $\eta \le S;$ so
from \eqref{meantaumu} we
also obtain:
\begin{equation} \label{eqcompatibility2}
 E( \tau ^ {(\mu )}_S) \le \frac 1 \mu E(S- \eta) .
\end{equation}
Letting $\mu$ go to zero in  \eqref{meantaumu}, we get:
\begin{equation} \label{meantaumu0}
E( \tau _S ^{(0 )} ) = - E( \eta ^2) + 2a E( \eta) + S(S-2a),
\end{equation}
and the compatibility condition for $\mu =0$  becomes
$$ - E( \eta ^2)+ 2a E(\eta) + S(S-2a) \ge 0 .$$
Finally, taking into account that $E(S- \eta ) \ge 0,$ we get from \eqref{meantaumu}:
\begin{equation} \label{eqcompatibility3}
 \mu E( \tau ^ {(\mu)} _S) \ge - \frac 1 { 2 \mu  } \ e^{ 2 \mu a } E \left ( e^{ - 2 \mu \eta } - e ^{ -2 \mu S } \right ).
\end{equation}
The  compatibility conditions above are
necessary so that the solution to the IFPT problem for regulated drifted BM exists, in the case of analytic Laplace transforms $\widehat f$ and $\widehat g.$
\begin{Remark}
{ \rm For regulated drifted BM $X^ {(\mu)},$  starting from $\eta \in [a,S],$
we have considered the FPT through $S$ from below.
In analogous way, if $\eta \in [S, b],$ one can consider $\widetilde \tau^ {(\mu )} _S ,$ namely  the
FPT of $X^ {(\mu)}$ through $S$ from above.
By using the same arguments of Proposition \ref{primoteorema}, and
\eqref{laplacetauabove},
it is possible to write the Laplace transform of $\widetilde \tau ^ {(\mu )} _S $ in terms of the Laplace transform
of  $\eta,$ if it exists; for instance, if $\mu =0,$
and we suppose that the density of $\eta $ is symmetric with respect to
$ (S+b)/2,$ then the solution $g  $ to the IFPT from above has the following Laplace transform (cf. with \eqref{laplacegmu0symm}):
\begin{equation}
\widehat g ( \theta ) = \frac {e^{-S \theta } + e ^{ -(2b - S) \theta } } {1 + e^{ (S-b) \theta } } \
\widehat f \left ( \frac {\theta ^2 } 2 \right ) , \ \theta \ge 0.
\end{equation}
}
\end{Remark}
\bigskip

Throughout the rest of the paper, for IFPT problem we will mean the problem concerning the FPT from below, namely,  $ \tau _S .$
\begin{Remark}
{ \rm Let $X(t)$ be regulated BM;  if $\tau_S$ has Gamma distribution,
there is no hope that there exists
the solution $g$  to the IFPT problem, with $g$
symmetric with respect to $(a + S)/2.$ In fact, suppose that
$f(t)$ is a Gamma density with some parameters $\alpha, \ \lambda >0,$ namely
$\widehat f ( \theta) = \left ( \frac \lambda {\theta + \lambda } \right ) ^ \alpha $ and
take $a=0$ and $S=1$ for the sake of simplicity;
then, inserting this expression of  $\widehat f ( \theta)$ into \eqref{laplacegmu0symm}, and calculating the second derivative at
zero of the candidate $\widehat g ,$ we see that it is $\le 0,$ implying that $\widehat g$ is not the Laplace transform
of any density $g$ of $\eta ,$ since it should be $E( \eta ^2)= \widehat g '' (0) \le 0,$ which is impossible.

}
\end{Remark}

Now, we further investigate the question of the existence of solutions to the IFPT problem. Referring to
regulated drifted BM,
we  will prove the existence of the density $g $ of the initial position $\eta \in [a,S]$
for a class of FPT densities $f.$
For the sake of simplicity, we limit ourselves to the case when $\mu =0, \ a=0, \ S=1< b $
and $g $ is required to be symmetric with respect to $1/2 ;$
in fact, for $\mu \neq 0$ the calculations involved are far more complicated. \par\noindent
For any integer $k \ge 0,$ set $I_k(\theta) = \int _{-1} ^1 e^ { - \theta x } x^k dx;$ as easily seen,
$I _0 (\theta )= 2 \sinh(\theta ) / \theta$ and the recursive relation
$ I _ {k} (\theta ) = \frac { (-1)^k e^ \theta - e ^ { - \theta } } \theta + \frac  k  \theta I _ {k-1} (\theta )$ allows to
calculate $I _k (\theta),$ for every $k.$ \par
The following Proposition gives a sufficient condition, in order that there exists the
solution to the  IFPT problem for regulated  BM.

\begin{Proposition} \label{existenceproposition}
Let $X$ be  regulated BM between the boundaries $0$ and $b,$ and let $S=1< b;$
suppose that the Laplace transform  of $f(t)$ has the form:
\begin{equation} \label{laplacedensityclass}
\widehat f (\theta) = \widehat f _ {2k}( \theta ) := \frac {\cosh (\sqrt { \theta/2} ) } {\cosh (\sqrt {2 \theta}) }
\left ( 1 + \frac 1 { 2k} \right ) \left [ \sqrt {\frac 2 \theta } \sinh \left ( \sqrt {\frac \theta 2 } \right ) -
I_{2k} \left ( \sqrt {\frac \theta 2 } \right ) \right ],
\end{equation}
for some integer $k>0.$
Then, there exists the solution $g=g_{2k}$ of the IFPT problem for $X,$  relative to the  FPT density $f,$ and
it results:
\begin{equation}
 g_{2k}(x)= \left (1 + \frac 1 { 2k} \right ) \left (1- (2x-1)^{2k} \right ), \ x \in (0,1).
\end{equation}
\end{Proposition}
{\it Proof.} We verify that $g_{2k}$ is the solution of the IFPT problem; a simple calculation shows that
$$\widehat g _ {2k} ( \theta) = \left (1 + \frac 1 {2k} \right ) e^{ - \theta /2}
\left [ \frac 2 \theta \sinh (\theta /2) - I_{2k} ( \theta /2) \right ] .$$
Since $g_{2k}$ is symmetric with respect to $1/2,$
the result follows by inserting $\widehat g _ {2k}$ into  \eqref{laplacegmu0symm}.
\par \hfill  $\Box$

\begin{Remark} \label{remdopoesistenza}
{ \rm A straightforward calculation shows that, if $\eta \in (0,1)$ has density $g_{2k},$ then $E(\eta ^2)= \frac {4k+5 } { 6(2k+3)} .$
By inserting $E( \eta ^2)$ into \eqref{meantaumu0} with $a=0, \ S=1,$  we obtain
that the FPT-distribution corresponding to $\widehat f _ {2k}$ has mean $E(\tau^ {(0)}_1 )= \frac {8k+13 } {6(2k+3) } \ .$
}
\end{Remark}

As an application of the results for regulated BM, we  consider now the piecewise-continuous process $\xi (t)$, obtained  by superimposing
to  BM a jump process, namely, for $\eta \in [a,S]$ and $t< T,$ we set $\xi (t)= \eta + B_t ,$
where $T$ is an exponentially distributed time with parameter $\lambda >0;$ we suppose that for $t=T$ the process $\xi (t)$ makes an
upward jump
and it crosses the barrier $S,$   irrespective
of its state before the occurrence of the jump. This kind of behavior is observed e.g.  in the presence of a so called {\it catastrophes} (see e.g.
\cite{dicre:queue03}).
Next, let us consider the reflected diffusion $\overline X$ with boundaries $a, \ b,$ associated to $\xi .$
Then, for $ \eta \in [a,S]$ the FPT of $\overline X$ over $S$
is
$\overline \tau _S = \inf \{ t>0: \overline X(t) \ge S \}.$
Conditionally on $\eta =x,$ we have:
$$P( \overline \tau _S(x) \le t )= P( \overline \tau _S (x)\le t | t<T ) P( t<T) + 1 \cdot P(t\ge T )=
P( \tau _S(x) \le t) e ^ {- \lambda t } + (1- e ^ {- \lambda t }). $$
Taking the derivative, we obtain  the FPT density of $\overline X,$  conditional to the starting position $x:$
$$ \overline f (t|x) = e ^ { - \lambda t } f(t|x) + \lambda e ^ { - \lambda t } \int _t ^ {+ \infty } f(s|x) ds .$$
By straightforward calculations, we obtain its Laplace transform:
$$ \widehat {\overline f } (\theta |x ) =  \int _0 ^ { \infty} e ^ { - \theta t } \overline f(t|x) dt=
\frac \theta { \lambda + \theta } \widehat f( \lambda + \theta |x ) + \frac  \lambda {\lambda + \theta } \ , \  \theta \ge 0 .$$
Finally, from the equation
$ \widehat {\overline f} (\theta ) = \int _a ^ S \widehat { \overline f } ( \theta |x) {\overline g}(x) dx,$  where
$\overline f$ is the FPT density of $\overline X,$ and $ \overline g$ is the density of $\eta, $
we get $\widehat {\overline f} .$ For the sake of simplicity, we limit ourselves to the case when $a= 0 $ and
the density $\overline g$ is symmetric with respect to $S/2 .$ Then, by using Theorem \ref{primoteorema}, we get:
$$ \widehat {\overline f} (\theta ) = \frac 1 { \lambda + \theta } \left [ \theta \cdot \frac {\widehat {\overline g}
\left ( \sqrt {2(\lambda + \theta) } \ \right) \left (1+ e^ {S\sqrt { 2(\lambda + \theta)}}  \ \right ) }   {2 \cosh \left ( S \sqrt {2(\lambda + \theta ) } \right )} + \lambda \right ] .$$
\par\noindent
Thus, by solving with respect to $ \widehat {\overline g} ( \theta),$ we  have obtained:
\begin{Proposition}
For $a=0 <S < b,$ if there exists a function  $\overline g ,$ symmetric with respect to $S/2 ,$
which is the solution to the  IFPT problem of $\overline X(t),$
relative to $S $ and
the FPT density $\overline f ,$
then its Laplace transform is given by:
\begin{equation} \label{laplacegconsalti}
 \widehat {\overline g } ( \theta )= \frac {2 \cosh (S \theta) }    {(\theta ^2 /2 - \lambda ) (1+ e ^{S \theta }) } \left [ \frac { \theta ^2 } 2 \
 \widehat {\overline f} \left ( \frac { \theta ^2 } 2 - \lambda \right ) - \lambda \right ] .
 \end{equation}
\end {Proposition}
\par
\hfill $\Box$

\begin{Remark}
{\rm For $\lambda =0,$ namely when no jump occurs, \eqref{laplacegconsalti} reduces to \eqref{laplacegmu0symm} with $a=0.$ }
\end{Remark}

\subsection{Reduction of reflected diffusions to reflected Brownian motion}
In certain cases, a reflected diffusion $ X$  can be reduced to  reflected BM by a variable change;
by using this approach,
we shall extend to reflected diffusions  the results obtained for reflected BM.
On the analogy of the definition holding for ordinary diffusions (see e.g. \cite{abundo:stapro12}, \cite{abundo:saa06}),
we introduce the following:

\begin{Definition}
Let $X(t)$ be a diffusion with reflecting boundaries $a$ and $b,$ which is driven by the SDER:
$$ dX(t)= \mu (X(t)) dt + \sigma (X(t)) dB_t  + dL_t - dU_t , \ X(0)=x \in [a,b] .$$
We say that $X(t)$
is conjugated to regulated BM if there exists
an increasing differentiable
function $V(x),$ with $V(0)=0,$
such that $X(t)=V^{-1} \left ( B_t + V(x) + \overline L _t - \overline U _t \right ),$ for any $t \ge 0,$
where $\overline L _t = V'(a) L_t$ and $\overline U _t = V'(b) U_t$ are regulators.
\end{Definition}
\bigskip

\noindent A class of reflected diffusions  which are conjugated to regulated BM is given by processes which are solutions
of SDERs such as:
\begin{equation} \label{conjdiffu}
 dX(t) = \frac 1 2 \sigma (X(t)) \sigma ' (X(t)) dt + \sigma (X(t)) dB_t + L_t - U_t, \ X(0)=x
\end{equation}
with $\sigma (\cdot) \ge 0.$
Indeed, if the integral
$V(x) :=  \int  ^ x \frac 1  {\sigma (r) }dr $
is convergent,
by It${\rm \hat o}$'s  formula for
reflected diffusions (see e.g. \cite{harr:85}), one gets  $V(X(t))= B_t +V(x) + V'(a) L_t - V'(b) U_t).$\bigskip

\noindent Let us consider a diffusion  $X,$
with reflecting boundaries $a$ and $b,$ which is conjugated to regulated BM via the function $V.$
Then, the process $Y(t):= V(X(t))$ is regulated BM between the boundaries $V(a)$ and $V(b),$
starting from $V(x),$ that is:
$$ Y(t)= V(x) + B_t + \overline L _t - \overline U _t \ ,$$
where $\overline L _t = V'(a) L_t$ and $\overline U _t = V'(b) U_t$ are the regulators of $Y(t),$ which increase only when
$Y=V(a)$ and $Y=V(b),$ respectively
Thus, for $x\in [a,S]:$
$$ \tau _S (x)= \inf \{ t \ge 0 : X(t) = S | X(0)=x \} = \inf \{ t \ge 0 : Y(t) = V(S) \}= \tau ^ Y_{S'} (V(x)),$$
where $S'= V(S)$ and  the superscript refers to the process $Y.$
Moreover,
if the initial position $\eta = X(0) \in [a,S] $  is random, its  density $g(x)$
is related to the corresponding density $\widetilde g $ of the initial position
$\widetilde \eta = V (\eta)\in [V(a), V(S)]$ of $Y(t),$  by the  equation
$\widetilde g  (y)=g(V^{-1} (y)) (V^{-1} )'(y),
\ y \in [V(a), V(S)].$ Furthermore,
the density $f(t)$ of $\tau _S$ and its Laplace transform are:
$$ f(t)= \int _{V(a)} ^{ V(S)} f ^Y (t | V(y)) \widetilde g  (y) dy  \ \ {\rm and} \ \ \
\widehat f(\theta)= \int _{V(a)} ^{ V(S)} \widehat f ^Y (\theta | V(y)) \widetilde g (y) dy ,$$
where $f ^Y (t | y)$ is the density of $ \tau ^ Y_{S'} (y)$ and $\widehat f ^Y (\theta | y)$ is its Laplace transform.
Therefore, if $X$ is conjugated to regulated BM via the function $V,$ then the solution $g$ to the IFPT problem for $X,$
relative to the FPT density $f$
and the barrier $S,$ can be written in terms of the solution $\widetilde g $ to the IFPT problem for regulated BM $Y(t)$
relative to the FPT density $f$
and the
barrier $V (S),$  by using that $g(x)= \widetilde g(V(x)) V'(x).$
From Theorem \ref{primoteorema} it follows that the Laplace transform of $\widetilde g$ satisfies
\eqref{laplacegmu0}, with $ \widehat g$ replaced by $\widehat {\widetilde g} , \ \eta$ replaced
with $\widetilde \eta , $  $a$ replaced with $V(a),$ and $S$ replaced with $V(S);$ in particular,
if one seeks  that $\widetilde g  $ is symmetric
with respect to $(V(a)+ V(S))/2,$ then (see \eqref{laplacegmu0symm})
the Laplace transform of $\widetilde g  $ is explicitly given by:
\begin{equation} \label{laplaceconj}
\widehat {\widetilde g }( \theta ) = \frac {e^{-V(S) \theta } + e ^{ -(2V(a) - V(S)) \theta } } {1 + e^{ (V(S)-V(a)) \theta } } \
\widehat f \left ( \frac {\theta ^2 } 2 \right ) , \ \theta \ge 0 .
\end{equation}
By inverting this Laplace transform, one can get $\widetilde g$ and therefore $g.$

\section{A few examples}
{\bf Example 1} Let $X(t)$ be regulated BM with boundaries $a, \ b \ (a <S<b),$ starting from $\eta \in [a, S]$ and
consider the FPT density:
\begin{equation} \label{densityfortheuniformex}
f(t) = \frac 1 {(S-a)^2} \sum _ {k=0} ^\infty \exp \left [  - \frac {\left (k+ \frac
1 2 \right )^2 \pi ^2 t} {2 (S-a)^2 } \right ],
\end{equation}
or the corresponding FPT Laplace transform:
\begin{equation} \label{laplaceexample1}
\widehat f ( \theta) = \frac {\tanh ( (S-a)\sqrt { 2 \theta} ) } {(S-a) \sqrt {2 \theta} } \ .
\end{equation}
Then, the solution $g$ to the IFPT problem for $X(t)$ is the uniform density in $(a,S).$
\par\noindent To verify this, we proceed backwardly; suppose that
$g(x)= {\bf 1} _ {(a,S)} (x) \frac 1 {S-a},$ so:
$$ \widehat g ( \theta ) = \int _ a ^S e^{ - \theta x } g(x) dx =
\frac {e^{ - a \theta  } - e ^{- S \theta } } { (S-a) \theta } \ .$$ Since $g$ is symmetric with respect $(a+S)/2,$
we can use \eqref{laplacegmu0symm}; after calculation, we obtain:
$$ \widehat f  (\theta ) = \frac {\sinh((S-a)\sqrt { 2 \theta} ) } {\sqrt {2 \theta} (S-a) \cosh( \sqrt {2 \theta}(S-a)) } \ ,$$
which yields \eqref{laplaceexample1}. In particular,
for $a= 0, \ S=1,$ \eqref{densityfortheuniformex} and \eqref{laplaceexample1} become:
\begin{equation} \label{example1a0}
f(t) = \sum _ {k=0} ^\infty \exp \left [  - \frac 1 2  \left (k+ \frac
1 2 \right )^2 \pi ^2 t   \right ] \ {\rm and } \
 \widehat f ( \theta) = \frac {\tanh ( \sqrt { 2 \theta} ) } { \sqrt {2 \theta} } \ .
 \end{equation}
The first three moments of this FPT distribution are $E \left [\tau ^{(0)} _1 \right ] = \frac 2 3 , \ E \left [ (\tau ^{(0)} _1)^2 \right ] = \frac {16 } {15 } , \
E \left [ (\tau ^{(0)} _1)^3 \right ] = \frac { 272} { 105 } ,$
and the solution $g$ to the IFPT problem is the uniform density in $(0,1).$\bigskip

\noindent {\bf Example 2}
For $a=0  < S < b,$ let $X(t)$ be regulated BM starting from $\eta \in [a,S],$ and
consider the FPT density whose Laplace transform is:
$$ \widehat f ( \theta) = \frac { \pi ^2} 2 \ \frac {1 + \cosh ( S \sqrt { 2 \theta } ) }
{\cosh ( S \sqrt { 2 \theta } ) (2 \theta S^2 + \pi ^2 ) } $$
(if e.g. $S=1,$ the Taylor expansion of $\widehat f$ up to the second order is
$ \widehat f ( \theta )= 1 - \left ( \frac 1 2 + \frac 2 {\pi ^2} \right ) \theta + \left ( \frac 5 {12} + \frac 1 {\pi ^2} +
\frac 4 {\pi ^4} \right ) \theta ^2 + o( \theta ^2), \ \theta \rightarrow 0 ).$
Then, the solution to the IFPT problem for $X$
is $g(x)= \frac { \pi } {2S} \sin \left ( \frac {\pi x} S  \right ), \ x \in (0,S) .$ \par\noindent
In fact, if we search for
a solution which is symmetric with respect to $S/2,$ we find by \eqref{laplacegmu0symm} that
$\widehat g (\theta ) = \frac {\pi ^2} 2 \ \frac {\left (1+e^{ - \theta S}\right ) } {\theta ^2 S^2 + \pi ^2 },$
which is indeed the Laplace transform of the function $g(x)$ above.
\bigskip

\noindent {\bf Example 3}
Take $a=0, \ S=1 ,$ and let $X(t)$ be regulated BM starting from $\eta \in [0,1];$
consider the FPT density whose Laplace transform is:
$$ \widehat f ( \theta) = \frac {\left (1 + e ^{ \sqrt { 2 \theta }} \right )
\left ( e ^{ \sqrt { 2 \theta}} -2 e ^{ \sqrt { \theta /2}} + 1 \right )   } {\theta \cosh ( \sqrt { 2 \theta} ) e^{ \sqrt { 2 \theta }} } $$
(the Taylor expansion of $\widehat f$ up to the second order is
$ \widehat f ( \theta )= 1 - \frac {17 } {24 }  \theta + \frac {811} {1440} \theta ^2 + o( \theta ^2), \ \theta \rightarrow 0 ).$ \par\noindent
Then, the solution to the IFPT problem for $X$
is  the triangular density in $[0,1]:$
$$g(x)=
\begin{cases}
4x, \ \ \ \ \ \ \ \ \ \ x \in [0, \frac 1 2 ] \\ 4x(1-x), \ x \in ( \frac 12 , 1 ]
\end{cases} \ .
$$
In fact, if we search for
a solution which is symmetric with respect to $1/2,$ we find by \eqref{laplacegmu0symm} that
$\widehat g (\theta ) = \frac 4 { \theta ^2 } \left ( 1 - e^{ - \theta /2 } \right ) ^2 ,$
which is indeed the Laplace transform of the function $g(x)$ above.
\bigskip

\noindent {\bf Example 4}
Take $a=0, \ S=1 ,$ and let $X(t)$ be regulated BM starting from $\eta \in [0,1];$
consider the FPT density whose Laplace transform is:
$$ \widehat f ( \theta) =
\frac {3 \left (1+ e ^{ \sqrt {2 \theta }}\  \right ) \left ( e^{- \sqrt {2 \theta } }
(\sqrt {2 \theta } +2) + \sqrt { 2 \theta} -2 \right ) }
{ \theta \sqrt {2 \theta } \left (e ^{ \sqrt { 2 \theta } } + e ^{ - \sqrt {2 \theta }} \ \right ) }  $$
(the Taylor expansion of $\widehat f$ up to the second order is
$ \widehat f ( \theta )= 1 - \frac {7 } {10 }  \theta + \frac {39} {70} \theta ^2 + o( \theta ^2), \ \theta \rightarrow 0 ).$ \par\noindent
Then, the solution $g$ to the IFPT problem for $X$
is  a Beta density in $[0,1],$ i.e. $ g(x)= 6 x(1-x).$
In fact, if we search for
a solution which is symmetric with respect to $1/2,$ we find by \eqref{laplacegmu0symm} that
$\widehat g (\theta ) = \frac 6 { \theta ^3 } \left ( (\theta +2) e^{ - \theta} + \theta - 2  \right )  ,$
which is indeed the Laplace transform of the function $g(x)$ above.
Notice that $\widehat f$ and $g$ are obtained as  special cases of $\widehat f_{2k}$ and $g_{2k}$ of
Proposition \ref{existenceproposition}, for $k=1.$
\bigskip

\noindent {\bf Example 5} For $0=a < S < b$ let $\overline X$ be the jump-process considered at the end of subsection 2.2., and
let:
$$ \widehat {\overline f } (\theta ) = \frac {1 } {\lambda + \theta } \left [
\frac { \theta \cdot \tanh ( S \sqrt {2(\lambda + \theta)}) } { S \sqrt {2(\lambda + \theta) } } + \lambda \right ]. $$
By Laplace inversion, one obtains:
$$ \overline f (t)= e^ { - \lambda t } \left [ \sum _ {k=0} ^ \infty \left ( \frac 1 {S^2} + \frac {2 \lambda S^2 } {(k+ \frac 1 2 ) ^2 \pi ^2  } \right )
\exp \left ( - \frac {(k+ 1/2 )^2 \pi ^2 t } {2 S^2} \right ) \right ],$$
which can be written as
\begin{equation}
\overline f (t)= e^ { - \lambda t } \phi (t) + \lambda e ^ { - \lambda t } \int _t ^ \infty \phi (s) ds ,
\end{equation}
where $\phi (t)$ is the FPT density considered in Example 1 with $a=0,$  i.e.:
\begin{equation}
 \phi (t)=  \frac 1 {S^2} \sum _ {k=0} ^ \infty \exp \left ( - \frac {(k+ 1/2 )^2 \pi ^2 t } {2 S^2  } \right ) .
 \end{equation}
Then, the solution $\overline g$  to the IFPT problem for $\overline X (t),$  relative to $S \in (0, b)$ and
 $\widehat {\overline f },$  which is  symmetric with respect to $S/2,$  is
the uniform density in $(0,S).$  \par\noindent
To verify this, it suffices to insert $\overline g (x) = {\bf 1} _ { (0,S) } (x) \frac 1 S$ into the equation
\eqref{laplacegconsalti} and to make some arrangements, by taking also into account the properties of Laplace transform.
\bigskip

\noindent {\bf Example 6} (Reflected geometric Brownian motion).
Let $0 < a < S < b,$ and $X(t)$  the solution of the SDER:
$$dX(t) = r X(t) dt + \sigma X(t) d B_t + dL_t - dU_t , \ X(0)= \eta \in [a,S]  , $$
where $ r$ and $\sigma$ are positive constant. The equation without reflecting is
well-known in the framework of Mathematical Finance,
since it describes the time evolution of a stock price. As easily
seen, $\ln X(t)= \ln \eta +  \mu t +  \sigma B_t +
\bar L _t - \bar U _t,$ where $\mu = r - \sigma ^2 /2 $ and $\bar
L _t, \bar U _t$ are  regulators; thus, $\ln X(t)/ \sigma$  is
regulated BM with drift $\mu / \sigma ,$ between the boundaries $\frac {\ln a } \sigma, \
\frac { \ln b } \sigma .$ Then, the IFPT problem for $X(t)$ relative to $S$
and the FPT density $f,$ is reduced to the IFPT problem for
regulated drifted BM, starting from $\frac {\ln \eta } \sigma ,$
relative to $\frac { \ln S }  \sigma$ and the same FPT density $f.$
Explicit examples can be obtained from Examples 1-4.
\bigskip

\noindent {\bf Example 7} Let $X$ be a reflected diffusion in
$[a,b],$ which is conjugated to regulated BM via the function $V;$
then, examples of solutions to the IFPT problem for $X$ can be
obtained from Examples 1-4 regarding regulated BM. \par\noindent
For instance, let us consider the FPT density
\begin{equation} \label{fptconjugated}
f(t) = \frac 1 {(V(S)-V(a))^2} \sum _ {k=0} ^\infty \exp \left [- \frac {  \left (k+ \frac
1 2 \right )^2 \pi ^2 t} {2 (V(S)-V(a))^2 } \right ],
\end{equation}
Then, the solution to the IFPT problem for $X$ relative to the
barrier $S \ (a <S< b) $ is:
\begin{equation} \label{gforconjugateuniform}
g(x)=  \frac { V'(x)}  {V(S)- V(a) } \cdot {\bf 1 } _ {(a,S)} (x).
\end{equation}
In fact, from \eqref{laplaceconj} and Example 1, $\widetilde \eta = V( \eta)$ turns out to be  uniformly distributed
in the interval $(V(a), V(S)),$
and so the relation $g(x)= \widetilde g(V(x)) V'(x)$ yields \eqref{gforconjugateuniform}.
\par
As explicit examples of reflected diffusions $X$ which are conjugated to regulated BM,
we mention the following.
\par\noindent
{\bf (i)} The process  driven by
\begin{equation}
 dX(t) = \frac 1 3 X(t)^ {1/3} dt + X(t)^ {2/3} \ dB_t  + dL_t - dU_t \ , \  X(0)=\eta \in [a,b],
\end{equation}
which is conjugated to regulated BM via the function $V(x) =3 x ^ {1/3}$
i.e. $X(t) = \left ( \eta ^ {1/3} + \frac 1 3 B_t  + \overline L_t - \overline U_t \right ) ^3.$
Here, as well as in the next examples,
$\overline L _t = V'(a) L_t$ and $\overline U _t = V'(b) U_t \ .$
\par\noindent
{\bf (ii)}  For $c>0,$ the process driven by
\begin{equation}
 dX(t) = \frac {3c^2} 8 ( X(t)  )^ {1/2} dt + c( X(t) ) ^ {3/4} \ dB_t +  dL_t - dU_t \ , \  X(0)=\eta \in [a,b] \ (a \ge 0),
\end{equation}
which is conjugated to regulated BM via the function $V(x) =\frac 4 c  x ^ {1/4}$ i.e.
$X(t) = \left ( \eta ^ {1/4} + \frac c 4 B_t  +  \overline L_t - \overline U_t \right ) ^4.$
\par\noindent
{\bf (iii)}  \ (Feller process or CIR model)\par\noindent
For $b> a \ge 0,$ the process  driven by
\begin{equation}
 dX(t) = \frac 1 4 dt + \sqrt {X(t) } \ dB_t +  dL_t - dU_t \ , X(0)=\eta \in [a,b] ,
\end{equation}
which is conjugated to regulated BM via the function
$V(x) =2 \sqrt x$ i.e. $X(t) = \frac 1 4 (B_t +2 \sqrt { \eta} +  \overline L_t - \overline U_t )^2.$ Notice that the process
is always $ \ge 0 .$
\par \noindent {\bf (iv)} (Wright \& Fisher-like process) \par\noindent
For $0 \le a < b \le 1,$ the process   driven by:
$$dX(t) = \left ( \frac 1 4 - \frac 1 2 X(t) \right ) dt + \sqrt {X(t)(1-X(t)) } \ dB_t   +  dL_t - dU_t\ , X(0) = \eta \in [a,b] ,  $$
which is conjugated to regulated BM via the function
$V(x) = 2 \arcsin \sqrt {x}.$
This equation is used for instance
in the Wright-Fisher model for population genetics and in certain
diffusion models for neural activity \cite{lansky:jtb94}; it results
$X(t) = \sin ^2 ( B_t /2 + \arcsin \sqrt {\eta} +  \overline L_t - \overline U_t)$ and so  $X(t)\in [0,1]$ for all $t\ge 0.$
Notice that, if we take $a=0$ and $b=1,$ both boundaries are attainable and there is no need for reflection in $a$ and $b,$
because the process without reflecting cannot exit the interval $[0,1],$ for any time $t$ (see e.g. \cite{abundo:stapro12}). \par
If the FPT density  is given by \eqref{fptconjugated},  from \eqref{gforconjugateuniform} we obtain
that the  solutions to
the IFPT problems  for the processes (i)--(iv) above, relative to the barrier $S,$ are explicitly given by: \par\noindent
$ g(x) = \left [3 x^{2/3} \left (S^{1/3} - a ^{1/3} \right ) \right ] ^{-1} \cdot {\bf 1 } _ {(a,S)} (x)$ (i), \
$ g(x) = \frac 1 4 \left [x^{3/4} \left ( S ^{1/4} -a ^{1/4} \right ) \right ] ^{-1} \cdot {\bf 1 } _ {(a,S)} (x)$ (ii), \par\noindent
$ g(x) = \frac 1 2 \left [\sqrt x \left ( \sqrt S  - \sqrt a  \right ) \right ] ^{-1} \cdot {\bf 1 } _ {(a,S)} (x)$ (iii), \
and
$ g(x) = \frac 1 2 \left [ \left ( \arcsin \sqrt S - \arcsin \sqrt a  \right ) \sqrt {x(1-x)} \  \right ] ^{-1} \cdot
{\bf 1 } _ {(a,S)} (x) $ (iv).

\end{document}